# A fast method for solving the linear optimization problem subjected to simplified Dombi-fuzzy relational equations

Amin Ghodousian , Sara Zal

*Abstract*—In this paper, an optimization model with a linear objective function subjected to a system of fuzzy relation equations (FRE) is studied where the feasible region is defined by the Dombi t-norm. Dombi family of t-norms includes a parametric family of continuous strict t-norms, whose members are increasing functions of the parameter. This family of t-norms covers the whole spectrum of t-norms when the parameter is changed from zero to infinity. Since the feasible solutions set of FREs is non-convex and the finding of all minimal solutions is an NP-hard problem, designing an efficient solution procedure for solving such problems is not a trivial job. Firstly, the feasible domain is characterized and then, based on some theoretical properties of the problem, a modified branch-and-bound solution technique is presented, which solves the problem by considering a few number of feasible paths. After presenting our solution procedure, a concrete example is included for illustration purpose.

## I. Introduction

In this paper, we study the following linear problem in which the constraints are formed as fuzzy relational equations defined by Dombi t-norm:

$$\begin{aligned} \min \quad & cx \\ & A\varphi x = b \\ & x \in [0,1]^n \end{aligned} \quad (1)$$

Where $I = \{1, \ldots, m\}$, $J = \{1, \ldots, n\}$, $A = (a_{ij})_{m \times n}$ is a fuzzy matrix such that $0 \leq a_{ij} \leq 1$ ($\forall i \in I$ and $\forall j \in J$), $b = (b_i)_{m \times 1}$ is a fuzzy vector such that $0 \leq b_i \leq 1$ ($\forall i \in I$), and "$\varphi$" is the Dombi t-norm defined as follows:

$$\varphi(x,y) = \begin{cases} 0 & x = 0 \text{ or } y = 0 \\ \left(1 + \left[\left(\frac{1-x}{x}\right)^\lambda + \left(\frac{1-y}{y}\right)^\lambda\right]^{\frac{1}{\lambda}}\right)^{-1} & \text{otherwise} \end{cases} \quad (2)$$

where $\lambda > 0$. If $a_i$ is the $i$'th row of matrix $A$, then problem (1) can be expressed as follows:

$$\begin{aligned} \min \quad & cx \\ & \varphi(a_i, x) = b_i \quad, \quad i \in I \\ & x \in [0,1]^n \end{aligned} \quad (3)$$

where the constraints mean

$$\varphi(a_i, x) = \max_{j \in J}\{\varphi(a_{ij}, x_j)\} = b_i \ (\forall i \in I) \text{ and}$$

$$\varphi(a_{ij}, x_j) = \begin{cases} 0 & a_{ij} = 0 \text{ or } x_j = 0 \\ \left(1 + \left[\left(\frac{1-a_{ij}}{a_{ij}}\right)^\lambda + \left(\frac{1-x_j}{x_j}\right)^\lambda\right]^{\frac{1}{\lambda}}\right)^{-1} & \text{otherwise} \end{cases} \quad (4)$$

As mentioned, the family $\varphi$ is increasing in $\lambda$. On the other hand, Dombi t-norm $\varphi$ converges to the basic fuzzy intersection $\min\{x, y\}$ as $\lambda$ goes to infinity and converges to Drastic product t-norm as $\lambda$ approaches zero. Therefore, Dombi t-norm covers the whole spectrum of t-norms [5]. In [3], the Dombi operations of single-valued neutrosophic numbers (SVNNs) were presented based on the operations of the Dombi t-norm and t-conorm. The authors proposed the single-valued neutrosophic Dombi weighted arithmetic average operator and the single-valued neutrosophic Dombi weighted geometric average operator to deal with the aggregation of SVNNs. In [1], a fuzzy morphological approach was presented to detect the edges of real time images in order to preserve their features, where Dombi t-

A. Gh. Faculty of Engineering Science, College of Engineering, University of Tehran, P.O.Box 11365-4563, Tehran, Iran(a.ghodousian@ut.ac.ir).

S. Z. Department of Engineering Science, College of Engineering, University of Tehran, Tehran, Iran(sara.zal@ut.ac.ir).
.

norm and t-conorm was used for computing morphological dilation and erosion. In [4], the authors studied the connection with Dombi aggregative operators, uninorms, strict t-norms and t-conorms. They presented a new representation theorem of strong negations that explicitly contains the neutral value.

The theory of fuzzy relational equations (FRE) was firstly proposed by Sanchez and applied in problems of the medical diagnosis [17]. Nowadays, it is well known that many issues associated with a body knowledge can be treated as FRE problems [16]. The solution set of FRE is often a non-convex set that is completely determined by one maximum solution and a finite number of minimal solutions [8]. The other bottleneck is concerned with detecting the minimal solutions that is an NP-hard problem [9,11,12,14]. The problem of optimization subject to FRE and FRI is one of the most interesting and on-going research topic among the problems related to FRE and FRI theory [2,6,8–12,15,18]. Recently, many interesting generalizations of the linear programming subject to a system of fuzzy relations have been introduced and developed based on composite operations used in FRE, fuzzy relations used in the definition of the constraints, some developments on the objective function of the problems and other ideas [7,11-13,15]. In this paper, after carefully study the solution set of system (1), we show that the problem can be solved by a modified branch-and-bound method.

The rest of the paper is arranged as follows. In Section 2, a special characterization of the feasible domain of problem (1) is derived. Based on the theoretical aspects of the problem, a new solutions set is obtained that includes all the minimal solutions. In Section 3, we study the effect of the cost vector c with the special characterization. A modified branch-and-bound method is presented and a step-by-step algorithm for solving problem (1) is given in Section 4. Some concluding remarks are made in Section 5.

## II. BASIC PROPERTIES OF MAX-DOMBI FRE Review Stage

Let $S$ denote the feasible solutions set of problem (1), that is, $S = \{x \in [0,1]^n : \max_{j=1}^{n}\{\varphi(a_{ij}, x_j)\} = b_i, \forall i \in I\}$. Also, for each $i \in I$, define $J_i = \{j \in J : a_{ij} \geq b_i\}$. Also, for each $i \in I$ and $j \in J_i$, we define

$$V(b_i, a_{ij}) = \left(1 + \left[\left(\frac{1-b_i}{b_i}\right)^\lambda - \left(\frac{1-a_{ij}}{a_{ij}}\right)^\lambda\right]^{\frac{1}{\lambda}}\right)^{-1} \quad (5)$$

According to [9], when $S \neq \emptyset$, it can be completely determined by one maximum solution and a finite number of minimum solutions. The maximum solution can be obtained by $\bar{X} = \min_{i \in I}\{\hat{x}_i\}$ where $\hat{x}_i = [(\hat{x}_i)_1, ..., (\hat{x}_i)_n]$ ($\forall i \in I$) is defined as follows

$$(\hat{x}_i)_j = \begin{cases} V(b_i, a_{ij}) & j \in J_i, \ b_i \neq 0 \\ 0 & j \in J_i, \ a_{ij} > b_i = 0 \\ 1 & \text{otherwise} \end{cases}, \ \forall j \in J \quad (6)$$

Moreover, if we denote the set of all minimum solutions by $\underline{S}$, then

$$S = \bigcup_{\underline{X} \in \underline{S}} \{x \in [0,1]^n : \underline{X} \leq x \leq \bar{X}\} \quad (7)$$

**Corollary 1.** $S \neq \emptyset$ if and only if $\bar{X} \in S$.

**Definition 1.** Let $i \in I$. For each $j \in J_i$, we define $\breve{x}_i(j) = [\breve{x}_i(j)_1, ..., \breve{x}_i(j)_n]$ such that

$$\breve{x}_i(j)_k = \begin{cases} V(b_i, a_{ik}) & b_i \neq 0 \text{ and } k = j \\ 0 & \text{otherwise} \end{cases}, \ \forall k \in J \quad (8)$$

**Definition 2.** Let $e: I \to \bigcup_{i \in I} J_i$ so that $e(i) = j \in J_i, \forall i \in I$, and let $E$ be the set of all vectors $e$. For the sake of convenience, we represent each $e \in E$ as an $m$-dimensional vector $e = [j_1, ..., j_m]$ in which $J_k = e(k)$.

**Definition 3.** Let $e \in E$. We define $\underline{X}(e) = [\underline{X}(e)_1, ..., \underline{X}(e)_n]$ where $\underline{X}(e)_j = \max_{i \in I}\{\breve{x}_i(e(i))_j\}, \forall j \in J$.

Based on the definitions 2 and 3, and according to [9], we have

$$\underline{S} \subseteq \{\underline{X}(e) : e \in E\} \quad (9)$$

Moreover, we have the following fundamental theorem [9]:

**Theorem 1.** Suppose that matrix $\bar{A} = (\bar{a}_{ij})_{m \times n}$ is resulted from $A = (a_{ij})_{m \times n}$ by the following steps:

(I) If $j_0 \notin J_i$ for some $i \in I$ and $j_0 \in J$, then we set $\bar{a}_{ij_0} = 0$.

(II) If $j_0 \in J_i$, $b_i \neq 0$ and there exists some $i' \in I$ ($i' \neq i$) such that $j_0 \in J_{i'}$, $b_{i'} \neq 0$ and $V(b_{i'}, a_{i'j_0}) < V(b_i, a_{ij_0})$, then we set $\bar{a}_{ij_0} = 0$.

(III) If $j_0 \in J_i$, $b_i \neq 0$ and there exists some $i' \in I$ ($i' \neq i$) such that $b_{i'} = 0$ and $a_{i'j_0} > 0$, then we set $\bar{a}_{ij_0} = 0$.
Then, two systems $S = \{x \in [0,1]^n : A\varphi x = b\}$ and $\bar{S} = \{x \in [0,1]^n : \bar{A}\varphi x = b\}$ have the same solutions.

**Definition 4.** For each $i \in I$, define $\bar{J}_i = \{j \in J : \bar{a}_{ij} \geq b_i\}$. Also, let $\bar{E}$ be the set of all vectors $e: I \to \bigcup_{i \in I} \bar{J}_i$ so that $e(i) = j \in \bar{J}_i, \forall i \in I$.

**Corollary 2.** $S \subseteq \{\underline{X}(e) : e \in \bar{E}\} \subseteq \{\underline{X}(e) : e \in E\}$.

**Proof.** Suppose that $\bar{a}_{ij_0} = 0$ by step (I). So, we have $j_0 \notin J_i$, i.e., $a_{ij_0} < b_i$. Thus, $\bar{a}_{ij_0} = 0 < b_i$ that means $j_0 \notin \bar{J}_i$. On the other hand, suppose that $\bar{a}_{ij_0} = 0$ by steps (II) or (III). So, $j_0 \in J_i$ ($a_{ij_0} \geq b_i$) and $b_i \neq 0$. In these cases, we have $\bar{a}_{ij_0} = 0 < b_i$, that is, $j_0 \notin \bar{J}_i$. Therefore, $\bar{J}_i \subseteq J_i$ that implies $\bar{E} \subseteq E$. Hence, $\{\underline{X}(e) : e \in \bar{E}\} \subseteq \{\underline{X}(e) : e \in E\}$. Now, the proof is resulted from relation (9) and Theorem 1. ☐

**Corollary 3.** $S = \bigcup_{e \in \bar{E}} [\underline{X}(e), \bar{X}]$.

**Proof.** The proof is resulted from relation (7) and Corollary 2. ☐

The following example illustrates the above-mentioned definitions.

**Example 1.** Consider the problem below with Dombi t-norm

$$\begin{bmatrix} 0.9452 & 0.4012 & 0.8976 & 0.6221 & 0.4368 & 0.8126 \\ 0.5271 & 0.1113 & 0.2456 & 0.3419 & 0.5271 & 0.2192 \\ 0.2073 & 0.8172 & 0.4386 & 0.4599 & 0.6152 & 0.2188 \\ 0.9111 & 0.7243 & 0.3274 & 0.8327 & 0.8327 & 0.5845 \end{bmatrix} \varphi \begin{bmatrix} 0.7243 \\ 0.5271 \\ 0.6152 \\ 0.8327 \end{bmatrix}$$

where

$$\varphi(x, y) = \begin{cases} 0 & x = 0 \text{ or } y = 0 \\ \left(1 + \sqrt{\left(\frac{1-x}{x}\right)^2 + \left(\frac{1-y}{y}\right)^2}\right)^{-1} & \text{otherwise} \end{cases}$$

(i.e., $\lambda = 2$). In this example, we have $J_1 = \{1,3,6\}$, $J_2 = \{1,5\}$, $J_3 = \{2,5\}$ and $J_4 = \{1,4,5\}$. According to relation (6), we attain $\hat{x}_1 = [0.7266\ 1\ 0.7335\ 1\ 1\ 0.7675]$, $\hat{x}_2 = [1\ 1\ 1\ 1\ 1\ 1]$, $\hat{x}_3 = [1\ 0.6312\ 1\ 1\ 1\ 1]$ and $\hat{x}_4 = [0.8506\ 1\ 1\ 1\ 1\ 1]$.

Therefore,
$\bar{X} = \min_{i \in I} \{\hat{x}_i\} = [0.7266\ 0.6312\ 0.7336\ 1\ 1\ 0.7675]$

Also, based on Definition 1, we have

$\check{x}_1(1) = [0.7266\ 0\ 0\ 0\ 0\ 0]$,
$\check{x}_1(3) = [0\ 0\ 0.7335\ 0\ 0\ 0]$,
$\check{x}_1(6) = [0\ 0\ 0\ 0\ 0\ 0.7675]$
$\check{x}_2(1) = [1\ 0\ 0\ 0\ 0\ 0], \check{x}_2(5) = [0\ 0\ 0\ 0\ 1\ 0]$
$\check{x}_3(2) = [0\ 0.6312\ 0\ 0\ 0\ 0]$,
$\check{x}_3(5) = [0\ 0\ 0\ 0\ 1\ 0]$
$\check{x}_4(1) = [0.8506\ 0\ 0\ 0\ 0\ 0]$,
$\check{x}_4(4) = [0\ 0\ 0\ 1\ 0\ 0], \check{x}_4(5) = [0\ 0\ 0\ 0\ 1\ 0]$

The cardinality of set $E$ is equal to $|E| = \prod_{i \in I} |J_i| = 3 \times 2 \times 2 \times 3 = 36$. So, we have 36 solutions $\underline{X}(e)$ associated to 36 vectors $e$. For example, for $e = [1,5,5,5]$, we obtain $\underline{X}(e) = \max\{\check{x}_1(1), \check{x}_2(5), \check{x}_3(5), \check{x}_4(5)\}$ from Definition 1,

2 and 3 that means $\underline{X}(e) = [0.7266\ 0\ 0\ 0\ 1\ 0]$. Now, from step (I) of Theorem 1, we have

$\bar{a}_{12} = \bar{a}_{13} = \bar{a}_{14} = \bar{a}_{15} = \bar{a}_{16} = \bar{a}_{22} = \bar{a}_{23} = \bar{a}_{24} = \bar{a}_{26} = \bar{a}_{31} = \bar{a}_{33} = \bar{a}_{34} = \bar{a}_{36} = \bar{a}_{42} = \bar{a}_{43} = \bar{a}_{46} = 0$

In all of these cases, $a_{ij} < b_i$, that is, $j \notin J_i$. Also, from part (II), we can set $\bar{a}_{21} = \bar{a}_{41} = 0$. For example, $a_{21} = b_2$ (i.e., $1 \in J_2$), $b_2 \neq 0$, $a_{11} > b_1$ (i.e., $1 \in J_1$), $b_1 \neq 0$ and $0.7266 = V(b_1, a_{11}) < V(b_2, a_{21}) = 1$.

Hence, we have

$$\bar{A} = \begin{bmatrix} 0.9452 & 0.0000 & 0.8976 & 0.0000 & 0.0000 & 0.8126 \\ 0.0000 & 0.0000 & 0.0000 & 0.0000 & 0.5271 & 0.0000 \\ 0.0000 & 0.8172 & 0.0000 & 0.0000 & 0.6152 & 0.0000 \\ 0.0000 & 0.0000 & 0.0000 & 0.8327 & 0.8327 & 0.0000 \end{bmatrix}$$

and therefore, $\bar{J}_1 = \{1,3,6\}$, $\bar{J}_2 = \{5\}$, $\bar{J}_3 = \{2,5\}$ and $\bar{J}_4 = \{4,5\}$ that imply $|\bar{E}| = \prod_{i \in I} |\bar{J}_i| = 3 \times 1 \times 2 \times 2 = 12$.

## III. LINEAR PROGRAMMING PROBLEM AND OPTIMAL SOLUTION

According to the well-known schemes used for optimization of linear problems such as (1) [5,10], problem (1) is converted to the following two sub-problems:

$$\min\ Z_1 = \sum_{j=1}^{n} c_j^+ x_j$$
$$\varphi(a_i, x) = b_i\ ,\ i \in I$$
$$x \in [0,1]^n$$

$$\min\ Z_2 = \sum_{j=1}^{n} c_j^- x_j$$
$$\varphi(a_i, x) = b_i\ ,\ i \in I$$
$$x \in [0,1]^n$$

where $c_j^+ = \max\{c_j, 0\}$ and $c_j^- = \min\{c_j, 0\}$ for $j = 1,2,\ldots,n$. It is easy to prove that $\bar{X}$ is the optimal solution of $Z_2$, and the optimal solution of $Z_1$ is $\underline{X}(e')$ for some $e' \in \bar{E}$.

**Theorem 2.** Suppose that $S \neq \emptyset$, and $\bar{X}$ and $\underline{X}(e^*)$ are the optimal solutions for $Z_2$ and $Z_1$, respectively. Then $c^T x^*$ is an optimal solution of problem (1), where $x^* = [x_1^*, x_2^*, \ldots, x_n^*]$ is defined as follows:

$$x_j^* = \begin{cases} \bar{X}_j & c_j \leq 0 \\ \underline{X}(e^*)_j & c_j > 0 \end{cases} \quad (10)$$

for $j = 1,2,\ldots,n$.

**Proof.** Let $x \in S$. Then, from Corollary 3 we have $x \in \bigcup_{e \in \bar{E}} [\underline{X}(e), \bar{X}]$. Therefore, for each $j \in J$ such that $c_j > 0$, inequality $x_j^* \leq x_j$ implies $c_j^+ x_j^* \leq c_j^+ x_j$. In addition, for each $j \in J$ such that $c_j \leq 0$, inequality $x_j^* \geq x_j$ implies $c_j^- x_j^* \leq c_j^- x_j$. Hence, $\sum_{j=1}^{n} c_j^- x_j^* \leq \sum_{j=1}^{n} c_j x_j$. □

As mentioned earlier, generating the maximum solution $\bar{X}$ is not a problem. If we know how to find a minimum solution for $Z_1$, then problem (1) can be solved via Theorem 2. In the next section, we provide a branch and bound method for solving $Z_1$.

## IV. MODIFIED BRANCH AND BOUND METHOD

A branch-and-bound method implicitly enumerates all possible solutions to a programming problem. For our application, in the beginning, we choose one constraint to branch the original problem into several sub-problems. Each sub-problem is represented by one node. Then branching at each node is done by adding one additional constraint. New sub-problems are created and represented by new nodes. Note that the more constraints added to a sub-problem, the smaller feasible domain it has and, consequently, the larger optimal objective value Z it achieves. Therefore, solving one sub-problem could provide the possibility to eliminate many possible solutions from consideration. In other words, once a current candidate solution is obtained, we judge other nodes for further consideration. If the best potential solution of one particular node is no better than the current candidate solution, then there is no need to branch on this node. Otherwise, branching is needed to yield a new bound. In this section, we use one concrete example to illustrate how this method works.

Based on the theory we built in previous sections, here we propose an algorithm for finding an optimal solution of problem (1).

**Step 1:** Compute $\bar{X} = \min_{i \in I} \{\hat{x}_i\}$.

**Step 2:** If $A\varphi\bar{X} \neq b$, according to Corollary 1, problem (1) has no feasible solution.

**Step 3:** Compute $\bar{J}_i, \forall j \in J$.

**Step 4:** Use the branch-and-bound concept based on $\{\underline{X}(e) : e \in \bar{E}\}$ to solve $Z_1$.

**Step 5:** Generate $x^*$ by (10).

**Example 2.** Consider problem (1) as follows:

min $\quad 6.2944 x_1 + 8.1158 x_2 - 7.4602 x_3 + 8.2675 x_4 + 2.6471 x_5 - 8.0491 x_6$

$\begin{bmatrix} 0.9452 & 0.0000 & 0.8976 & 0.0000 & 0.0000 & 0.8126 \\ 0.0000 & 0.0000 & 0.0000 & 0.0000 & 0.5271 & 0.0000 \\ 0.0000 & 0.8172 & 0.0000 & 0.0000 & 0.6152 & 0.0000 \\ 0.0000 & 0.0000 & 0.0000 & 0.8327 & 0.8327 & 0.0000 \end{bmatrix} \varphi x = \begin{bmatrix} 0.7243 \\ 0.5271 \\ 0.6152 \\ 0.8327 \end{bmatrix}$

$x \in [0,1]^6$

**Step 1:** We find
$\bar{X} = [0.7266 \ 0.6312 \ 0.7336 \ 1 \ 1 \ 0.7675]$.

**Step 2:** Since $A\varphi\bar{X} = b$, we know $S \neq \emptyset$.

**Step 3:** As mentioned before, we have $J_1 = \{1,3,6\}$, $J_2 = \{1,5\}$, $J_3 = \{2,5\}$ and $J_4 = \{1,4,5\}$ and also, $\bar{J}_1 = \{1,3,6\}, \bar{J}_2 = \{5\}, \bar{J}_3 = \{2,5\}$ and $\bar{J}_4 = \{4,5\}$. Hence, the cardinality of $\{\underline{X}(e) : e \in E\}$ and $\{\underline{X}(e) : e \in \bar{E}\}$ is equal to $|E| = \prod_{i \in I} |J_i| = 36$ and $|\bar{E}| = \prod_{i \in I} |\bar{J}_i| = 12$, respectively.

**Step 4:** Sub-problem $Z_1$ is obtained as follows:

min $\quad Z_1 = 6.2944 \, x_1 + 8.1158 \, x_2 + 8.2675 \, x_4 + 2.6471 \, x_5$

$\begin{bmatrix} 0.9452 & 0.0000 & 0.8976 & 0.0000 & 0.0000 & 0.8126 \\ 0.0000 & 0.0000 & 0.0000 & 0.0000 & 0.5271 & 0.0000 \\ 0.0000 & 0.8172 & 0.0000 & 0.0000 & 0.6152 & 0.0000 \\ 0.0000 & 0.0000 & 0.0000 & 0.8327 & 0.8327 & 0.0000 \end{bmatrix} \varphi x = \begin{bmatrix} 0.7243 \\ 0.5271 \\ 0.6152 \\ 0.8327 \end{bmatrix}$

$x \in [0,1]^6$

Our branch-and-bound method begins by choosing **j** from $\bar{J}_1 = \{1,3,6\}$. Hence, each feasible solution must satisfy either $e(1) = 1$, $e(1) = 3$ or $e(1) = 6$. This yields three branches denoted by nodes 1, 2 and 3 in Fig. 1. We then obtain a lower bound on the Z-value associated with each node. For example, when $e(1) = 1$, we know $c_1 V(b_1, a_{11}) = 6.2944(0.7266) = 4.5735$, therefore, any choice with $e(1) = 1$ results in $Z \geq 4.5735$. So we write $Z \geq 4.5735$ for node 1 of Fig. 1. Similarly, any choice with $e(1) = 3$ or $e(1) = 6$ results in $Z \geq 0$. Since at this moment we have no reason to exclude any one of nodes 1, 2 and 3 from consideration, all the nodes have to be further investigated. By using the jump-tracking technique, we branch on the node with a lower bound on Z. In this ease, we choose node 2. Since $\bar{J}_2 = \{5\}$, any choice associated with node 2 must satisfy $e(2) = 5$. Branching on node 2 yields node 4 in Fig. 1. For this new node, we need to evaluate a lower bound for the objective value Z.
So, at node 4, we compute

$Z \geq c_3 V(b_1, a_{13}) + c_5 V(b_2, a_{25})$
$\quad = 0 \, (0.7266) + 2.6471(1) = 2.6471$

Hence, any choice associated with node 4 must have $Z \geq 2.6471$. Again, we do not have any reason to exclude node 4 from consideration, so we need to branch on one node. So, we branch node 4. Any choice associated with node 4 must satisfy either $e(3) = 2$ or $e(3) = 5$, since $\bar{J}_3 = \{2,5\}$.

This yields nodes 5 and 6 in Fig. 1. By the same reasoning as before, any choice associated with node 5 must have

$Z \geq c_3 V(b_1, a_{13}) + c_5 V(b_2, a_{25}) + c_2 V(b_3, a_{32}) = 0\,(0.7266) + 2.6471(1) + 8.1158(0.6312) = 7.7697$

and any choice associated with node 6 must have

$Z \geq c_3 V(b_1, a_{13}) + c_5 \max\{V(b_2, a_{25}) + V(b_3, a_{35})\} = 0\,(0.7266) + 2.6471 \max\{1,1\} = 2.6471$

.
Of course, we are interested in node 6. To branch further on node 6, any choice associated with node 6 must have $e(4) = 4$ or $e(4) = 5$, since $\bar{J}_4 = \{4,5\}$. This yields nodes 7 and 8 in Fig. 1. Note that node 7 corresponds to the sequence $e = [3,5,5,4]$. This sequence leads to a value of Z= $10.9146$. Therefore, node 7 is a feasible sequence which may be viewed as a candidate solution with Z= $10.9146$.

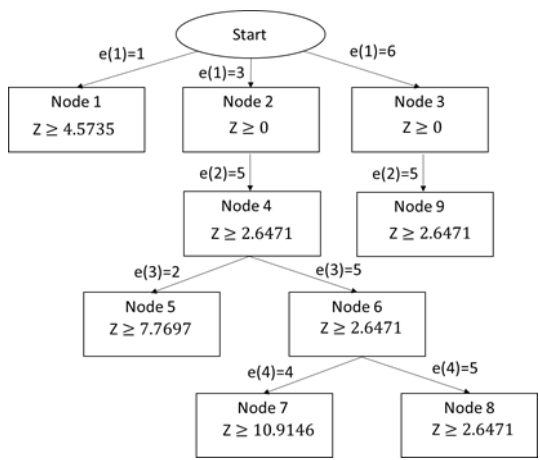

Figure 1. Modified branch and bound method.

By the same manner, node 8 corresponds to the sequence $e = [3,5,5,5]$ leading to a value of Z= $2.6471$. Because the Z value for nodes 1 and 7 cannot be lower than 2.6471, these nodes can be eliminated from further consideration. However, node 3 cannot be eliminated yet, because it is still possible for node 3 to yield a sequence having Z< $2.6471$. Hence, we now branch on node 3. Since $\bar{J}_2 = \{5\}$, any sequence associated with node 3 must have $e(2) = 5$. Correspondingly, we create node 9. For node 9, we calculate Z≥ $2.6471$. Since any sequence associated with node 9 must have Z≥ $2.6471$, this node is eliminated from consideration. Now, with the exception of node 9, every node in Fig. 1 has been eliminated from consideration. Node 8 yields the sequence $e = [3,5,5,5]$ that generates solution $\underline{X}(e) = [0\ 0\ 0.7335\ 0\ 1\ 0]$. Thus, we should choose the sequence $e = [3,5,5,5]$ and solution

$\underline{X}(e) = [0\ 0\ 0.7335\ 0\ 1\ 0]$ with Z= $2.6471$ as an optimal solution for $Z_1$.

**Step 5:** According to (10), we can find the optimal value of the original problem as $x^* = [0\ 0\ 0.7336\ 0\ 1\ 0.7675]$ with optimal objective value $Z^* = -9.0033$.

V. CONCLUSIONS

In this paper, we have studied a linear optimization problem subjected to a system of Dombi fuzzy relation equations and presented a procedure to find an optimal solution. After analyzing the properties of its feasible domain, we found a solutions set that include all the minimal solutions of the main feasible region, then applied a modified branch-and-bound method to find one solution. From the analysis of Theorem 2, it is clearly seen that if all minimum solutions of a given system of fuzzy relation equations can be found, then an optimal solution of the optimization problem defined by (10) can be constructed. Extension to other types of objective functions is under investigation.

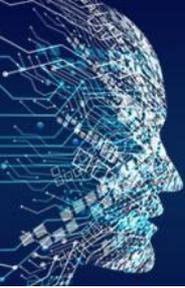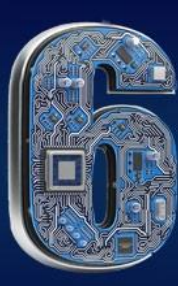